\documentclass[12pt,a4paper]{amsart}

\usepackage[left=3cm,right=3cm,top=3.5cm,bottom=3.5cm]{geometry}

\usepackage{amsmath,amsthm,amssymb}
\usepackage[english]{babel}

\usepackage{tikz,pgffor}

\newtheorem{theorem}{Theorem}[section]
\newtheorem{lemma}[theorem]{Lemma}
\newtheorem{prop}[theorem]{Proposition}

\newtheorem{coro}[theorem]{Corollary}
\newtheorem{conj}[theorem]{Conjecture}
\theoremstyle{definition}
\newtheorem{dfn}[theorem]{Definition}
\newtheorem{rem}[theorem]{Remark}
\newtheorem{exam}[theorem]{Example}

\numberwithin{equation}{section}

\DeclareMathOperator*{\Spec}{Spec}

\DeclareMathOperator{\rk}{rk}

\newcommand{\rleft}{\mathopen{}\mathclose\bgroup\left}
\newcommand{\rright}{\aftergroup\egroup\right}

\def\T{{\Bbb T}}
\def\C{{\Bbb C}}

\def\Q{{\Bbb Q}}
\def\R{{\Bbb R}}
\def\P{{\Bbb P}}
\def\Z{{\Bbb Z}}
\def\A{{\Bbb A}}

\begin{document}

\title[Projecting lattice polytopes according to MMP]
{Projecting lattice polytopes according to the Minimal Model Program}

\author{Victor V. Batyrev}
\address{Fachbereich Mathematik, Universit\"at T\"ubingen, Auf der
Morgenstelle 10, 72076 T\"ubingen, Germany}
\curraddr{}
\email{batyrev@math.uni-tuebingen.de}
\thanks{}

\begin{abstract}
The Fine interior $F(P)$ of a $d$-dimensional 
lattice polytope 
$P \subset \R^d$ is the set of all points  $y \in P$ having 
integral distance at least $1$ to any integral supporting 
hyperplane of $P$.  We call a lattice polytope $F$-hollow if its Fine interior is empty. The main theorem claims that up to unimodular equivalence in each dimension $d$ there exist  only finitely many  $d$-dimensional $F$-hollow    
lattice polytopes $P$, so called {\em sporadic}, which do not admit   a lattice projection 
onto a $k$-dimensional $F$-hollow  lattice polytope $P'$ for some $1 \leq k \leq d-1$. The proof is purely combinatorial, but  it is inspired by $\Q$-Fano fibrations in the Minimal Model Program, since we show that  
non-degenerate toric hypersurfaces $Z \subset (\C^*)^d$ defined by zeros of Laurent polynomials with a given Newton polytope $P$ have  negative Kodaira dimension if and only if $P$ is $F$-hollow.  The finiteness theorem for 
$d$-dimensional sporadic $F$-hollow Newton polytopes $P$ give rise to 
finitely many families ${\mathcal F}(P)$ of $(d-1)$-dimensional $\Q$-Fano hypersurfaces with at worst canonical  singularities.
\end{abstract}
\maketitle

\thispagestyle{empty}

\section{Introduction}

Let $P \subset \R^d$ be an arbitrary 
$d$-dimensional convex polytope.

\begin{dfn}
For a nonzero lattice vector ${\bf a}=(a_1, \ldots, a_d) \in \Z^d \setminus 
\{0 \}$ consider  
 the {\bf integral supporting hyperplane} of $P$: 
 $$H_{{\bf a}, P} := \left\{ {\bf x}= (x_1, \ldots, x_d) \in \R^d\, \mid \, \sum_{i=1}^d a_i x_i = {\rm Min}_P({\bf a}) :=\min_{{\bf y} \in P} \left( \sum_{i=1}^n a_i y_i \right)  \right\}.$$ 
If ${\rm gcd}(a_1, \ldots, a_d) =1$  we call the number 
\[{\rm dist}_\Z({\bf y}, H_{{\bf a}, P}) :=  \sum_{i=1}^n a_i y_i  \; - \; {\rm Min}_P({\bf a})  \geq 0 \]
the {\bf integral distance} between   
a   point ${\bf y} =(y_1, \ldots, y_d) \in P$  and  integral  supporting hyperplane $H_{{\bf a},P}$.  
\end{dfn} 

\begin{dfn} \cite{Fine83,Reid87}
Let $P \subset \R^d$ be an arbitrary 
$d$-dimensional convex polytope. The 
 set $F(P)$ of all points in $P$ having integral distance at least $1$ to any integral supporting hyperplane $H_{{\bf a}, P}$ is called the {\bf Fine interior of $P$}, i.e., 
 \[ F(P) := \{ {\bf y} \in \R^d \, \mid \, {\rm dist}_\Z ({\bf y}, H_{{\bf a}, P}) \geq 1, \;\; \forall {\bf a}  \in \Z^d \setminus \{0\} \}. \]  
\end{dfn}

\begin{rem}
If the affine span of every facet of $P$ is a integral supporting 
hyperplane, then  Gordan's lemma  shows that 
among  countably many inequalities ${\rm dist}_\Z ({\bf y}, H_{{\bf a}, P}) \geq 1$  defining $F(P)$ only finitely many ${\bf a} \in \Z^d$ 
are necessary. In particular,  
$F(P) \subset P$ is a rational polytope (or empty set) 
if all  vertices of $P$ belong to $\Q^d$.  We explain more 
details concerning this fact in Proposition \ref{toric1}. {\bf Lattice polytopes $P$}, i.e., polytopes 
having vertices in $\Z^d$, are main objects of our study.  However, for some technical 
reasons it will be convenient to consider the Fine interior $F(P)$ 
of rational polytopes $P$ and the Fine interior $F(\lambda P)$ of their arbitrary  real positive multiples $\lambda P$ $(\lambda \in \R_{>0}$). 
\end{rem}

\begin{rem} \label{int-point}
If $P  \subset \R^d$ is a $d$-dimensional lattice polytope, then every interior lattice point 
${\bf m} \in Int(P) \cap \Z^d$ necessarily belongs to $F(P)$, 
and we obtain  the inclusion 
\begin{equation*}
{\rm Conv}\,(Int(P) \cap \Z^d) \subset F(P)
\end{equation*} 
which  is in fact equality for lattice polytopes $P$ of dimension $d \in \{1, 2\}$ \cite{Ba23}.
\end{rem}

Recall some standard definitions. 

\begin{dfn}
Two lattice polytopes
$P_1, P_2 \subset \R^d$ are called 
{\bf unimodular equi\-valent} if there exists
a lattice-preserving affine isomorphism $\varphi\, :\, 
\R^d \to \R^d$ 
such that $\varphi(P_1) = P_2$.
\end{dfn} 

\begin{dfn} A $k$-dimensional lattice 
polytope $P' \subset \R^{k}$ $(1 \leq k  < d)$ is called  {\bf lattice projection}, or {\bf $\Z$-projection},   
of a $d$-dimensional lattice polytope $P \subset \R^d$ if there 
exists an affine map $\pi \, :\, \R^d \to \R^{k}$ inducing   
a surjective  map of  lattices  $\pi \, :\, \Z^d \to \Z^{k}$ and $\pi(P) = P'$. 
\end{dfn}

\begin{rem} 
Since every epimorphism  
$\pi \, :\, \Z^d \to \Z^{k}$ splits,  we can choose a splitting 
$\Z^d \cong \Z^{k} \oplus \Z^{n-k}$ such that the $\Z$-projection $\pi$ has 
the standard form:
\[ \pi (x_1, \ldots, x_d) = (x_1, \ldots, x_{k}) \in \R^{k} \; \; \forall 
(x_1, \ldots, x_d) \in \R^d. \]
\end{rem}

\begin{dfn}
We call a $d$-dimensional polytope $P \subset \R^d$ 
{\bf $F$-hollow} if $F(P) = \emptyset$. 
\end{dfn}

\begin{rem}
Let   $P' \subset \R^k$ be a lattice projection 
of a $d$-dimensional polytope $P \subset \R^d$, then 
the lattice epimorphism $\pi \, :\, \Z^d \to \Z^{k}$ allows 
to lift every integral supporting hyperplane of $P'$ to 
an integral supporting hyperplane of $P$. In particular, 
the condition $F(P') = \emptyset$ implies $F(P) = \emptyset$. This fact  allows to construct  infinitely many pairwise unimodular distinct  $d$-dimensional $F$-hollow 
lattice polytopes $P$ whose lattice projections are equal to a given lower-dimensional $F$-hollow lattice polytope $P'$. 
\end{rem}

\begin{dfn}
We call a $d$-dimensional $F$-hollow lattice polytope $P$ 
{\bf sporadic}, if $P$  does not admit any lattice projection $\pi\, : P \to P'$ 
onto a $k$-dimensional $F$-hollow lattice polytope $P'$ $(1 \leq k \leq d-1)$. 
\end{dfn}

The present paper shows that in any fixed dimension $d$, 
apart from  finitely many unimodular equivalence 
classes of $d$-dimensional sporadic  $F$-hollow lattice polytopes,  
every $d$-dimensional  $F$-hollow lattice polytope $P$ admits  a $\Z$-projection to some 
lower-dimensional   $F$-hollow lattice  
polytope $P'$:  

\begin{theorem} \label{proj-F}
In each dimension $d$ there exist up to 
unimodular transformations  
only finitely many $d$-dimensional sporadic $F$-hollow lattice polytopes $P$. 
\end{theorem}

\begin{center}
\begin{tikzpicture}[scale=0.8, xshift= -5.0cm]

\draw[fill=blue!30!white, very thick] (2,0) -- (3,0) --(2,1)--(2,0);

\draw[fill=blue!30!white, very thick] (-3,0) -- (1,0) --(-3,1)--(-3,0);

\draw[fill=blue!30!white, very thick] 
(-2,-1) -- (1,-1) -- (3,-2) -- (-2,-2) -- (-2,-1);

\draw[step=1cm,gray,very thin]
(-4.4,-2.4) grid (4.4,2.4) ;

\fill (2,0) circle (3pt); 
\fill (2,1) circle (3pt); 

\fill (0,-1) circle (3pt);
\fill (-1,-1) circle (3pt);
\fill (-2,-1) circle (3pt);
\fill (1,-1) circle (3pt);

\fill (-2,0) circle (3pt);

\fill (3,0) circle (3pt);  %

\fill (-3,1) circle (3pt);
\fill (-3,0) circle (3pt);

\fill (1,0) circle (3pt);
\fill (0,0) circle (3pt); 
\fill (-1,0) circle (3pt);

\fill (1,-2) circle (3pt);
\fill (0,-2) circle (3pt); 
\fill (-1,-2) circle (3pt);
\fill (2,-2) circle (3pt);
\fill (3,-2) circle (3pt); 
\fill (-2,-2) circle (3pt);

\node [left] at (4.5,-3) {{ countably many $F$-hollow  polygons}}; 
\begin{scope}[xshift= 8cm] 
\draw[fill=blue!30!white, very thick] (1,-1) -- (-1,1) --(-1,-1)--(1,-1);
\fill (0,0) circle (3pt);
\fill (-1,1) circle (3pt);
\fill (-1,0) circle (3pt);
\fill (-1,-1) circle (3pt); 
\fill (0,-1) circle (3pt);
\fill (1,-1) circle (3pt);

\node [left] at (2.5,-3) {{a sporadic polygon}}; 
\draw[step=1cm,gray,very thin]
(-2.4,-2.4) grid (2.4,2.4) ;
\end{scope}
\end{tikzpicture}
\end{center}

\begin{rem} Theorem 
 \ref{proj-F} is a purely combinatorial contribution to the theory of lattice polytopes but its motivation  comes 
 from birational algebraic geometry \cite{Ba23,Tre10}. 
Recall that a lattice $P$ is called {\bf hollow} if 
$P$ contains no lattice points in its interior. 
By \ref{int-point},  any $F$-hollow lattice polytope $P$ is hollow. In dimension $d \in \{1, 2 \}$ the inverse statement is also true.  Therefore,  for $d \leq  3$, 
Theorem \ref{proj-F} follows from 
a result of Treutlein 
\cite{Tre10} which  was later  generalized  by Nill and Ziegler for hollow lattice polytopes in arbitrary dimension $d$ \cite{NZ11}.  
In case  $d \geq 4$,   
Theorem \ref{proj-F} does not follow from  Theorem of Nill and Ziegler  because hollow lattice polytopes $P$  of  dimension $\geq 3$  need not be $F$-hollow. For instance, there exist up to unimodular equivalence exactly $9$ examples of  $3$-dimensional hollow lattice polytopes which are not $F$-hollow \cite[Appendix B]{BKS22}.  
\end{rem}

\section{The proof of Theorem \ref{proj-F}}

The proof of Theorem  \ref{proj-F} uses standard  combinatorial notions from the theory 
of toric varieties \cite{CLS11}. 
We consider two dual to each other 
lattices $M \cong \Z^d$ and 
$N:= Hom(M, \Z)$ together with their scalar 
extensions $M_\R:= M \otimes \R$, $N_\R := N \otimes \R$
and the natural pairing 
\[  \langle * , * \rangle \; : \; M_\R \times N_\R \to \R. \]
Every $d$-dimensional rational polytope $P$ defines a  convex piecewise 
linear function  
\[ {\rm Min}_P \, : \, N_\R \to \R, \;\;  y \mapsto \min_{u \in P} \langle u, y \rangle \]
whose domains of linearity form   a complete 
rational polyhedral 
fan $\Sigma_P$, a collection   
of rational polyhedral cones $\sigma$  in $N_\R$, which is called the 
{\bf normal fan of $P$}. 

Using the above notations, the  Fine interior $F(P)$ of 
a  polytope $P \subset M_\R$ 
can be equivalently reformulated  as follows:

\begin{dfn}
\[ F(P) := \left\{ x \in M_\R \, \mid \, \langle x, \nu \rangle  \geq {\rm Min}_P(\nu) \, + 1,  \,  \forall \nu \in N \setminus \{ 0 \} \right\} . \]
\end{dfn}

Our first idea in the proof of \ref{proj-F} is to consider a  positive number $\mu = \mu(P)$ attached to $P$: 

\begin{dfn}
Let $P \subset M_\R$ be a $d$-dimensional rational polytope. 
We call the number  
\[ \mu(P):= {\rm inf}\, \{ \lambda \in \R_{>0} \; \mid \; 
F(\lambda P) \neq \emptyset \} \]
the {\bf minimal multiplier of $P$}.     
\end{dfn} 

\begin{rem}
If $P \subset M_\R$ is a $d$-dimensional lattice polytope, then it follows from standard properties of Ehrhart polynomials that the lattice polytope $(d+1)P$ always contains at least one interior lattice point 
and hence $F((d+1)P) \neq \emptyset$. This implies the 
 inequality $$\mu(P) \leq d +1 $$ 
for all  $d$-dimensional lattice polytopes $P$. 
Note that the   inequality 
is sharp, because $\mu(P) = \dim P +1$ if 
$P$ is the $d$-dimensional lattice simplex 
spanned 
by $0 \in \Z^d$ and by the standard lattice basis $e_1, \ldots, e_d$ of $\Z^d$. 
\end{rem} 

\begin{rem} It is clear that  $\mu(P) \leq  1$ if and only if $P$ is not $F$-hollow. 
\end{rem}

We will use the following property of the minimal multiplier $\mu(P)$ which will 
be proved later in Propositions \ref{toric2} and \ref{fano-comb}:

\begin{prop} \label{mult}
Let $P \subset M_\R$ be a $d$-dimensional rational polytope. Then the following statements hold:

{\rm (i)} The number $\mu(P)$ is rational. 

{\rm (ii)} For a positive rational number $\lambda$, one has 
$\lambda = \mu(P)$ if and only if 
$$0 \leq \dim F(\lambda P) \leq d -1.$$ 

{\rm (iii)} If  $\dim F(\mu P) =0$, then the convex 
hull $Q:= {\rm Conv}(S_F(\mu P))$ of the set 
\[ S_F(\mu P):= \{ \nu \in N \; \mid {\rm Min}_{F(\mu P)}(\nu) = {\rm Min}_{\mu P}(\nu) +1 \}  \subset N \]
is a $d$-dimensional lattice polytope in $N_\R$ having 
a unique interior lattice point $0 \in N$.  
\end{prop}

Finally, we need  Theorem  of Hensley \cite{Hen83} and some its generalizations due to  Lagarias and Ziegler \cite{LZ91}.

\begin{theorem} \cite{Hen83} \label{th-Hen83}
For any given positive integers $k, d$, there exists a constant $C(k, d)$ depending only on $k$ and $d$ such that the volume 
$vol(P)$ of any $d$-dimensional lattice polytope having 
exactly $k$ interior lattice points is bounded from above 
by $C(k,d)$.  
\end{theorem}

\begin{theorem}\cite[Thm. 2]{LZ91} \label{th-LZ91}
If a $d$-dimensional lattice polytope $P$ has volume $vol(P) \leq V$, then $P$  is unimodular equivalent to a lattice polytope in the 
$d$-dimensional lattice cube 
\[ \{ {\bf x} \in \R^d \; \mid \; 
0 \leq x_i \leq n \cdot n! V, \; i=1, \ldots, d \}. \] 
In particular, a family of $d$-dimensional lattice 
polytopes $P_i$ $(i \in I)$ contains only finitely many unimodular equivalence classes if and only if there exists a constant 
$C>0$ such that  
\[ vol(P_i) < C \;\; \forall i \in I.\] 
\end{theorem}

\noindent
{\bf Proof of Theorem \ref{proj-F}.} Let $P \subset M_\R$ be a   $d$-dimensional $F$-hollow lattice polytope. 
We put   $\mu := \mu(P)$.  
By \ref{mult}(ii),  we have $0 \leq \dim F(\mu P) \leq d-1$. 
Consider two cases (two pictures below illustrate the case $d =2$).

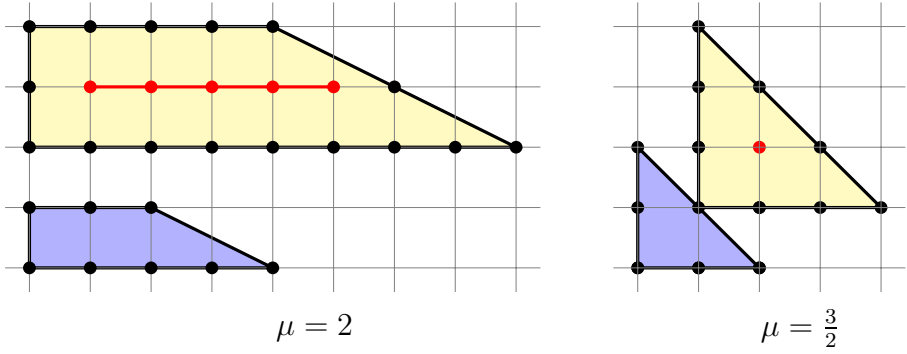
\begin{center}
\begin{tikzpicture}[scale=0.8]
\draw[fill=blue!30!white, very thick] 
(-4,-2) -- (0,-2) -- (-2,-1) -- (-4,-1) -- (-4,-2);

\draw[fill=yellow!30!white, very thick] 
(-4,0) -- (-4,2) -- (0,2) -- (4,0) -- (-4,0);

\draw[step=1cm,gray,very thin]
(-4.4,-2.4) grid (4.4,2.4) ;

\draw[red, very thick] 
(-3,1) -- (1,1);

\fill (0,2) circle (3pt); 
\fill (-1,2) circle (3pt); 
\fill (-2,2) circle (3pt); 
\fill (-3,2) circle (3pt); 
\fill (-4, 2) circle (3pt); 

\fill (-4,0) circle (3pt); 
\fill (4,0) circle (3pt); 

\fill (2,0) circle (3pt); 
\fill (2,1) circle (3pt);

\fill (-2,0) circle (3pt);

\fill (3,0) circle (3pt);  %

\fill (-4,1) circle (3pt);
\fill (-3,0) circle (3pt);

\fill (1,0) circle (3pt);
\fill (0,0) circle (3pt); 
\fill (-1,0) circle (3pt);

\fill (-4,-1) circle (3pt);
\fill (-3,-1) circle (3pt);
\fill (-2,-1) circle (3pt); 

\fill (0,-2) circle (3pt);
\fill (-1,-2) circle (3pt);
\fill (-2,-2) circle (3pt);
\fill (-3,-2) circle (3pt); 
\fill (-4,-2) circle (3pt);

\fill[color=red] (1,1) circle (3pt);
\fill[color=red] (0,1) circle (3pt);
\fill[color=red] (-1,1) circle (3pt);
\fill[color=red] (-2,1) circle (3pt);
\fill[color=red] (-3,1) circle (3pt);

\node [left] at (1.5,-3) {{\bf $\mu =2$}}; 
\begin{scope}[xshift= 8cm] 
\draw[fill=yellow!30!white, very thick] (-1,2) -- (2,-1) --(-1,-1)--(-1,2);
\draw[fill=blue!30!white, very thick] (0,-2) -- (-2,0) --(-2,-2)--(0,-2);
\fill (-1,-1) circle (3pt);
\fill (-2,0) circle (3pt);
\fill (-2,-1) circle (3pt);
\fill (-2,-2) circle (3pt); 
\fill (-1,-2) circle (3pt);
\fill (0,-2) circle (3pt);

\fill[color=red] (0,0) circle (3pt);

\fill (-1,0) circle (3pt);
\fill (0,-1) circle (3pt);
\fill (-2,-2) circle (3pt); 
\fill (-1,-2) circle (3pt);
\fill (0,-2) circle (3pt);
\fill (2,-1) circle (3pt);
\fill (0,-1) circle (3pt);
\fill (1,0) circle (3pt);
\fill (-1,-0) circle (3pt); 
\fill (-1,2) circle (3pt);
\fill (0,1) circle (3pt);
\fill (-1,1) circle (3pt);
\fill (1,-1) circle (3pt);

\node [left] at (1.5,-3) {{$\mu = \frac{3}{2}$}};  
\draw[step=1cm,gray,very thin]
(-2.4,-2.4) grid (2.4,2.4) ;
\end{scope}
\end{tikzpicture}
\end{center}

{\bf Case 1. $ 1 \leq \dim F(\mu P) \leq d-1$.}
We define the sublattice $N' \subset N$ consisting of 
all $\nu \in N$ such that $\langle x, \nu \rangle = \langle x', \nu \rangle $ for any two points $x, x' \in F(\mu P)$. Then 
$N/N'$ has no torsion elements and $N'$, because 
if $l\nu \in N'$ for some positive integer $l$, then   $\langle x, l \nu \rangle = \langle x', l \nu \rangle $ for any two points $x, x' \in F(\mu P)$, hence $\langle x, \nu \rangle = \langle x', \nu \rangle $, i.e., $\nu \in N'$. Therefore  $N'$ is 
is a direct summand of $N$, and    the embedding $N' \hookrightarrow N$ defines a lattice projection 
$\pi\,:\, M \to M':= Hom(N', \Z), $  
where
 \[ 1 \leq {\rk}(M') =
d - \dim F(\mu P) \leq  d-1. \] 
Consider the lattice polytope $P':= \pi(P)$. It remains to show that $F(P') =\emptyset$. By definition of $N'$,  
$\pi (F(\mu P))$ is some rational point $q \in M'_\Q$. Moreover, one has 
$q= F(\mu \pi(P)) =F(\mu P')$. Since $\mu >1$, 
by monotonicity of Fine interior \cite[Remark 3.7]{Ba23}, the  polytope $F(P')$ must be strictly smaller than the point $q =F(\mu P')$. Hence   $F(P') = \emptyset$, i.e., $P'$ is $F$-hollow. In fact we have shown that $\mu(P') = \mu(P) = \mu$
and $\dim F(\mu P')= 0$. 
\medskip

{\bf Case 2. $\dim F(\mu P) = 0$.} By \ref{mult}(i),  
 $F(\mu P)$ is a rational point  $p \in M_\Q$
 and for any $\nu \in S_F(\mu P)$ we have 
 \[ 
 {\rm Min}_{F(\mu P)}(\nu) =
  \langle p, \nu \rangle = 1 +
  {\rm Min}_{\mu P}(\nu) = 1 + \min_{x \in \mu 
  P} \langle x, \nu \rangle. \]
Equivalently,  we have 
$ \langle x, \nu \rangle  \geq -1 + \langle p, \nu \rangle$ for all 
 $x \in \mu P$ and for all  $\nu \in S_F(\mu P),$ 
or 
\[  \langle x, \nu \rangle  \geq -1  \;\; 
\forall x \in \mu P - p , \, 
\forall \nu \in S_F(\mu P).  \]
The last conditions imply the inequalities 
\[ \langle x, y \rangle  \geq -1   \;\; 
\forall x \in \mu P - p , \, 
\forall y  \in Q := {\rm Conv}( S_F(\mu P)).   \] 
and we obtain the inclusion $\mu P - p \subseteq Q^*$, 
where 
\[ Q^* := \{ x \in M_\R\, \mid \, \langle x, y \rangle \geq -1 , \; \forall 
y \in Q \} \]
is the rational polar dual polytope of  $Q$. 
By \ref{mult}(iii),  
$Q \subset N_\R$ is a $d$-dimensional lattice polytope  
having only one interior lattice point $0$. 
By theorems of Hensley \ref{th-Hen83} and Lagarias-Ziegler \ref{th-LZ91}, up to unimodular transformations 
there exist only finitely many 
possibilities for the $d$-dimensional lattice polytope 
$Q$. Therefore, 
there exists only finitely many possible values 
for volumes $vol(Q^*)$ of the dual polytope of $Q$, i.e., $vol(Q^*)$ is bounded
by some constant $C(d)$ depending only on $d$.   Since $\mu >1$, 
it follows from 
the inclusion $\mu P - p \subseteq  Q^*$ that 
$$C(d) \geq vol(Q^*) \geq vol(\mu P -p) = vol(\mu P)  = 
\mu^d vol(P)> vol(P).$$ 
Hence, by Theorem \ref{th-LZ91} of Lagarias and Ziegler, up to unimodular equivalence we obtain only finitely many of $d$-dimensional 
lattice polytopes $P$ such that  $\dim F(\mu P)=0$. 
\hfill $\Box$

\begin{dfn}
We call a $d$-dimensional $F$-hollow lattice polytope $P \subset \R^d$  {\bf weakly sporadic}, if $\dim F(\mu P) = 0$, where $\mu$ is the minimal multiplier $\mu(P)$ of $P$. 
\end{dfn}

\begin{coro}
In each  dimension $d$   there exist up to unimodular equivalence only finitely many 
weakly sporadic $d$-dimensional $F$-hollow lattice polytopes. Moreover, if a $F$-hollow lattice polytope $P$ with minimal multiplier $\mu$ is not sporadic, then $P$ admits a canonical lattice projection $\pi\,:\, P \to P'$, 
where $P'$ is a $k$-dimensional weakly sporadic 
$F$-hollow lattice polytope $(1 \leq k \leq d-1)$ with 
the same minimal multiplier $\mu = \mu(P) = \mu(P')$.  
\end{coro} 
 
 \noindent {\em Proof.}
The statemen immediately follows from Cases 1 and  2 in the proof of Theorem \ref{proj-F}. \hfill $\Box$  
 
\begin{rem}
  If $d \geq 2$, one can easily find examples of $d$-dimensional weakly sporadic $F$-hollow  lattice polytopes $P$ which are not sporadic  
$F$-hollow polytopes, i.e. $P$ admitting lattice projections onto a $k$-dimensional 
$F$-hollow polytope $P'$ $(1 \leq k \leq d-1)$. 
\end{rem}

\begin{exam}
Up to unimodular equivalence there exist exactly three  weakly sporadic $F$-hollow lattice polygons $P$  which are not sporadic $F$-hollow lattice polygons: 

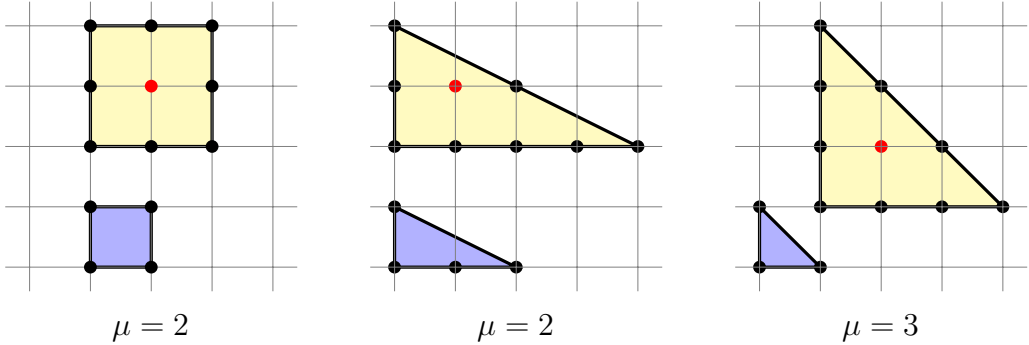
\begin{center}
\begin{tikzpicture}[scale=0.8]
\draw[fill=blue!30!white, very thick] 
(0,-1) -- (0,-2) -- (-1,-2) -- (-1,-1) -- (0,-1);
\draw[fill=yellow!30!white, very thick] 
(-1,0) -- (1,0) -- (1,2) -- (-1,2) -- (-1,0);

\draw[step=1cm,gray,very thin]
(-2.4,-2.4) grid (2.4,2.4) ;

\fill (0,2) circle (3pt); 
\fill (0,0) circle (3pt); 
\fill (-1,2) circle (3pt); 
\fill (1,2) circle (3pt); 
\fill (-1, 0) circle (3pt); 
\fill (1,0) circle (3pt);
\fill (1,1) circle (3pt); 
\fill (-1,1) circle (3pt); 

\fill[color=red] (0,1) circle (3pt);

\fill (0,-1) circle (3pt); 
\fill (0,-2) circle (3pt); 
\fill (-1,-2) circle (3pt); 
\fill (-1,-1) circle (3pt); 

\node [left] at (0.8,-3) {{\bf $\mu =2$}}; 
\begin{scope}[xshift= 6cm] 
\draw[fill=yellow!30!white, very thick] (-2,0) -- (2,0) --(-2,2)--(-2,0);
\draw[fill=blue!30!white, very thick] (-2,-1)--(0,-2) -- (-2,-2) --(-2,-1);

\fill (-2,0) circle (3pt);
\fill (-1,0) circle (3pt);
\fill (0,0) circle (3pt);
\fill (1,0) circle (3pt);
\fill (2,0) circle (3pt);

\fill[color=red] (-1,1) circle (3pt);
\fill (-2,1) circle (3pt);
\fill (0,1) circle (3pt);

\fill (-2,2) circle (3pt);

\fill (-2,-1) circle (3pt);
\fill (-2,-2) circle (3pt);
\fill (-1,-2) circle (3pt);
\fill (0,-2) circle (3pt);

\node [left] at (0.8,-3) {{$\mu = 2$}}; 
\draw[step=1cm,gray,very thin]
(-2.4,-2.4) grid (2.4,2.4) ;
\end{scope}
\begin{scope}[xshift= 12cm] 
\draw[fill=yellow!30!white, very thick] (-1,2) -- (2,-1) --(-1,-1)--(-1,2);
\draw[fill=blue!30!white, very thick] (-2,-1)--(-1,-2) -- (-2,-2) --(-2,-1);

\fill (0,-1) circle (3pt);
\fill (-1,0) circle (3pt);
\fill (-2,-1) circle (3pt);
\fill (1,0) circle (3pt);
\fill (2,-1) circle (3pt);

\fill[color=red] (0,0) circle (3pt);
\fill (1,-1) circle (3pt);
\fill (0,1) circle (3pt);

\fill (-1,2) circle (3pt);

\fill (-1,-1) circle (3pt);
\fill (-1, 1) circle (3pt);
\fill (-1,-2) circle (3pt);
\fill (-2,-2) circle (3pt);

\node [left] at (0.8,-3) {{$\mu = 3$}}; 
\draw[step=1cm,gray,very thin]
(-2.4,-2.4) grid (2.4,2.4) ;
\end{scope}
\end{tikzpicture}
\end{center}
\end{exam}

\begin{exam} \label{Del-Pezzo-123}
It follows from the combinatorial classification of all maximal hollow $3$-dimensional lattice polytopes obtained by Averkov et. al. \cite{AWW11,AKS17} 
that the following three $3$-dimensional weakly sporadic   $F$-hollow lattice  with the minimal multipliers $\mu \in 
\{  \frac{7}{6},  \frac{5}{4}, \frac{4}{3} \}$:
\[ \Delta_1 := {\rm Conv}\{ (0,0,0), (2,0,0), (0,3,0), (0,0,6) \}, \; \mu = \frac{7}{6}; \]
 \[ \Delta_2 := {\rm Conv}\{ (0,0,0), (2,0,0), (0,4,0), (0,0,4) \},  \; \mu = \frac{5}{4};  \]
 \[ \Delta_3 := {\rm Conv}\{ (0,0,0), (3,0,0), (0,3,0), (0,0,3) \},  \; \mu = \frac{4}{3}  \]
are in fact sporadic $F$-hollow polytopes. 
\end{exam}

\begin{rem}
We will see in the last section that every $d$-dimensional 
 weakly sporadic  $F$-hollow polytope 
$P$  defines a family ${\mathcal F}(P)$ of  non-degenerate $(d-1)$-dimensional $\Q$-Fano toric hypersurfaces 
with at worst canonical singularities. For example,  three lattice tetrahedra  $\Delta_i$, $i \in \{ 1,2,3\}$  from Example \ref{Del-Pezzo-123} are Newton polytopes  of smooth Del Pezzo surfaces of the anticanonical degree $i \in   \{ 1,2,3\}$ naturally imbedded into $3$-dimensional toric weighted projective 
spaces $\P(1,1,2,3)$, $\P(1,1,1,2)$, and $\P(1,1,1,1)$ respectively. We  expect that the complete list of all unimodular classes 
of  $3$-dimensional weakly sporadic  $F$-hollow lattice 
polytopes $P$  must have reasonable length.  
We draw attention to the fact that this list   includes not only   $\Delta_1, \Delta_2, \Delta_3$ 
 but also $31$ more $3$-dimensional weakly sporadic    $F$-hollow polytopes $P$  
with $\mu(P) =2$ arising  from  $3$-dimensional Gorenstein 
polytopes of index $2$ which  were classified in \cite{BJ10}.  
\end{rem}

\section{The minimal multiplier $\mu(P)$}

Let $P \subset M_\R$ be a $d$-dimensional rational  polytope. 
Theory of toric varieties associates with the normal fan 
$\Sigma_P$ a $d$-dimensional 
projective toric variety $X_P$ together  
with the ample $\Q$-Cartier divisor 
\[ L_P := \sum_{\nu \in \Sigma_P[1]} -{\rm Min}_P(\nu) D_\nu, \]
where $D_\nu$ $(\nu \in \Sigma_P[1])$ are torus invariant divisors on $X_P$ 
corresponding to primitive lattice generators  $\nu$ of
$1$-dimensional cones in $\Sigma_P$,  and we have  
\[ P = \{ x \in M_\Q \; \mid \; \langle x, \nu \rangle \geq 
{\rm Min}_P(\nu), \; \forall \nu \in \Sigma_P[1]\} . \] 

Theory of toric varieties allows to describe the Fine 
interior $F(P)$ as rational polytope associated 
with the  adjoint divisor on some smooth projective toric variety $Y_\Sigma$ obtained by a regular 
refinement $\Sigma$ of the normal fan $\Sigma_P$. More precisely, one has 

\begin{prop}\label{toric1}
Let $\rho\, :\, Y= Y_\Sigma \to X_P$ be a projective desingularization of the toric variety $X_P$ corresponding 
to a regular simplicial 
refinement $\Sigma$ of the normal fan $\Sigma_P$. Denote   by $D_{\nu_1}, \ldots, D_{\nu_s}$ the  torus invariant divisors corresponding 
to primitive lattice vectors in $\Sigma[1]= \{ \nu_1, \ldots, \nu_s \}$. Then the Fine interior 
$F(P)$ is the rational polytope defined by $s$ inequalities 
\[ \langle x , \nu_i \rangle \geq {\rm Min}_P(\nu_i) +1, \;\; i \in \{1, \ldots, s\} \]
corresponding to the adjoint 
divisor on $Y$:
\[ K_Y + \rho^*(L_P) =  \sum_{\nu  \in \Sigma[1]} \left(-1 - {\rm Min}_{P} (\nu) \right) D_{\nu}. \] 
\end{prop}

\noindent
{\em Proof.} Let $\nu \in N$ an arbitrary nonzero lattice vector. Then there exists a minimal regular simplicial cone $\sigma \in \Sigma$ containing $\nu$. Without loss of generality, we can assume that $\nu_1, \ldots, \nu_r \in  \Sigma[1]$ $(r \leq d <s)$ are generators of $\sigma$. 
Then $\nu = \sum_{i=1}^r l_i \nu_i$ for some positive integer
coefficients $l_1, \ldots, l_r$. Since $\Sigma$ is a refinement of $\Sigma_P$, there exists a minimal $r'$-dimensional cone $\sigma' \in \Sigma_P$  containing $\sigma$ ($r \leq r'\leq d)$. The cone $\sigma'$ is dual to some $(d-r')$-dimensional face $\Theta' \prec P$ and we have  
\[ {\rm Min}_P(\nu_i) = \langle x, \nu_i \rangle, \;\; \forall x \in \Theta', \;  \forall i =1, \ldots, r. \] 
Since ${\rm Min}_P( \cdot)$ is linear  on $\sigma \subset \sigma'$, we obtain  
\[ {\rm Min}_P(\nu) = {\rm Min}_P\left(\sum_{i=1}^r l_i \nu_i \right) =
\sum_{i=1}^r l_i {\rm Min}_P(\nu_i),   \]
 and $r$ inequalities appearing in definition of   $F(P)$
 \[ \langle x , \nu_i \rangle \geq {\rm Min}_P(\nu_i) +1, \;\; i \in \{1, \ldots, r\}  \]
imply that 
 \[   \langle x, \nu \rangle  =   
 \sum_{i=1}^r l_i  \langle x, \nu_i \rangle \geq 
 \sum_{i=1}^r l_i{\rm Min}_P(\nu_i) + \sum_{i=1}^r l_i  = {\rm Min}_P(\nu) + \sum_{i=1}^r l_i \geq {\rm Min}_P(\nu) +1,  \]
and the equality $\langle x, \nu \rangle = {\rm Min}_P(\nu) +1$ for some $x \in F(P)$ can happen only if $r=1$ and $\nu = \nu_1 \in \Sigma[1]$.
Therefore, $s$  inequalities 
\[ \langle x , \nu_i \rangle \geq {\rm Min}_P(\nu_i) +1, \;\; i \in \{1, \ldots, s\} \]
are already sufficient to obtain the Fine interior $F(P)$. 
Finally, we note 
that the  canonical divisor of toric variety $Y= Y_\Sigma$ equals $- \sum_{i=1}^s D_{\nu_i}$ and $\Q$-Cartier divisor $\rho^*L_P$ equals $-\sum_{i=1}^s {\rm Min}_P(\nu_i)D_{\nu_i}$. Therefore,  
the rational polytope $F(P)$ 
corresponds to the adjoint divisor $K_Y + \rho^*L_P$ of $Y$.
\hfill $\Box$

\begin{prop} \label{toric2}
Let $X_P$ be projective toric variety corresponding to a $d$-dimensional rational polytope $P$. Consider any projective toric  desingularization
$\rho\, :\, Y= Y_\Sigma \to X_P$    as in \ref{toric1}.  Denote by  $\Lambda_{\rm eff}(Y) 
\subset {\rm Pic}(Y)_\R$ 
the closed cone of effective divisors  of the smooth projective toric 
variety $Y$. Let  
$L:= \rho^*L_P \in {\rm Pic}(Y)_\Q$ be the pullback of the ample $\Q$-Cartier 
divisor $L_P$.  Then 
\[ \mu(P) =  {\rm inf}\, \{ \lambda \in \R_{>0} \, \mid \, 
 [K_Y] + \lambda [L] \in \Lambda_{\rm eff}(Y)\},  \]
is a rational number and for $\mu:= \mu(P)$ one has 
\[ 0 \leq \dim F(\mu P) < d. \] 
\end{prop}

\noindent 
{\em Proof.}  Using Cox coordinates on $Y$ one easily obtains that the cone of effective divisors
$\Lambda_ {\rm eff}(Y) \subset {\rm Pic}(Y)_\R$ is a  
rational polyhedral cone generated  by the classes $[D_\nu]$ $(\nu \in \Sigma[1])$
of torus invariant divisors.  The class $[L]$ represenst the class of  a semiample 
big $\Q$-Cartier divisor  $\rho^*(L_P)$ on $Y$ which defines a 
rational point $[\rho^*(L_P)] \in {\rm Pic}(Y)_\R$
in the interior of the cone $\Lambda_ {\rm eff}(Y)$. 
On the other hand, the canonical class $[K_Y] \in {\rm Pic}(Y)_\R$ of toric variety 
$Y$ does not belong to the cone $C_{\rm eff}(Y)$, because $Y$ is a rational toric variety.  Therefore,  the ray 
\[ \{ [K_Y] + \lambda [\rho^* L_P] \, \mid \, \lambda \in \R_{\geq 0} \}  \subset {\rm Pic}(Y)_\R \] 
with the rational origin $[K_Y]$ having the  rational direction $[\rho^* L_P]$ must hit the rational polyhedral cone $\Lambda_ {\rm eff}(Y)$ 
in some rational point $[K_Y] + \mu[\rho^* L_P] \in \Lambda_ {\rm eff}(Y)$ located at the polyhedral boundary 
$\partial \Lambda_ {\rm eff}(Y)$ contained in some 
proper rational polyhedral  facet $\Gamma \prec \Lambda_ {\rm eff}(X)$. Since the intersection point $[K_Y] + \mu [\rho^* L_P]$
of the rational ray $[K_Y] + \lambda [\rho^* L_P]$ with the facet $\Gamma$ has rational coordinates, the number 
$\mu$ must be rational. Now we use the fact that the class $[D]$ of a $\Q$-divisor
\[ D= \sum_{i=1} b_i [D_{\nu_i}] \in {\rm Pic}(Y)_\Q, \;\; b_i \in \Q,  \]
represents a point in $\Lambda_{\rm eff}(Y)$ if and only if 
the rational polytope 
\[   P_D: = \{ x \in M_\R\, \mid \, \langle x, \nu_i \rangle \geq - b_i, \;\; \forall \, i \in \{1, \ldots, s \} \} \]
is not empty. Moreover, $[D]$ represents an interior point in 
$\Lambda_{\rm eff}(Y)$ if and only if the rational polytope $P_D$ has maximal dimension $d$. Applying \ref{toric1} to the adjoint $\Q$-divisor 
$$D:= K_Y + \mu \rho^*L_P = K_Y + \rho^*\mu L_P = 
K_Y + \rho^* L_{\mu P},$$
we obtain that $F(\mu P)$ is not empty and $\dim F(\mu P) <d$.
\hfill  $\Box$

\begin{dfn}
Let $P \subset M_\R$ be an arbitrary  $d$-dimensional rational 
polytope. Assume $F(P) \neq \emptyset$. Then we call the 
set 
\[ S_F(P):= \{ \nu \in N  \, \mid \, 
{\rm Min}_{F(P)}(\nu) 
=  {\rm Min}_P(\nu) +1 \} \]
the {\bf support of the Fine interior} of $P$. It follows from \ref{toric1} that $S_F(P)$ is always a finite set whose positive convex span $\R_{\geq 0} S_F(P)$ equals $N_\R$
since 
\[ F(P) = \{ x \in M_\R\, \mid \, \langle x, \nu \rangle \geq {\rm  Min}_P(\nu) +1 \;\; \forall \nu \in S_F(P) \}\]
is compact. 
\end{dfn}

\begin{prop} \label{fano-comb}
Let $P$ be a  $d$-dimensional rational polytope with  $\dim F(P) =0$. Then the convex 
hull $Q:= {\rm conv}(S_F(P)) \subset N_\R$ is a $d$-dimensional lattice 
polytope containing in its interior only one lattice point $0 \in N$.
\end{prop}

\noindent
{\em Proof.}  Let  $p:= F(P) \in M_\Q$. Note  that the shifted polytope $P_0:= P-p$ has the Fine interior $0 \in M$, and the sets 
$S_F(P)$ and $S_F(P_0)$ are the same, 
since $F(P_0) = F(P) -p$. Hence we can assume $F(P) = p = \{0\} \subset Int(P)$. The lattice polytope $Q ={\rm conv}(S_F(P))$ is $d$-dimensional, since   
\[ 0 = F(P) = \{ x \in M_\R\, \mid \, \langle x, \nu \rangle \geq 0 \;\; \forall \nu \in S_F(P) \}. \]
 Using the upper 
convex piecewise linear function ${\rm Min}_P\, :\, N_\R \to \R$, we obtain the dual to $P$ rational polytope  
\[ P^*:= \{ y \in N_\R \, \mid \, {\rm Min}_P(y) \geq  -1\}. \] 
Since  ${\rm Min}_{F(P)}(\nu) = \langle 0, \nu \rangle =0$ for all $\nu \in N$, we obtain $S_F(P) = \partial P^* \cap N$, where $\partial P^* := \{ y \in N_\R \, \mid \, {\rm Min}_P(y) =  -1\} $ is the boundary of $P^*$. The $d$-dimensional rational polytope $P^*$ has only $0 \in N$  as interior lattice point, because ${\rm  Min}_P(\nu) \leq -1$ for all $\nu \in N \setminus \{0\}$. Hence the $d$-dimensional lattice subpolytope $Q \subseteq P^*$ has also  only $0 \in N$ as its interior lattice point.  
\hfill  $\Box$

\begin{rem}
We note that the minimal multiplier $\mu(P)$ has naturally  
appeared in the  arithmetical problem of counting rational points of bounded height on algebraic varieties \cite{BaMa90}. The close  relation between  the boundary 
point $[K_Y + \mu L] \in \Lambda_{\rm eff}(Y)$ and the Minimal Model Program was observed  in \cite{Ba92}. Fano fibrations of  smooth toric varieties $Y$ associated with adjoint divisors 
$K_Y + \mu L$ were considered in \cite{BTsch96,BTsch98}. 
\end{rem}

\section{Toric hypersurfaces with Newton polytope $P$}

Let $A \subset M \cong  \Z^d$ be a finite subset such that 
the convex hull $P:= {\rm Conv}(A) \subset M_\R \cong \R^d$ is a 
$d$-dimensional lattice polytope. Take an arbitrary field $K$ and consider $P$ 
as {\bf Newton polytope} of a Laurent 
polynomial 
$$f ({\bf t}) = \sum_{{\bf m} \in A} {c}_{\bf m} {\bf t}^{\bf m} \in 
K[t_1^{\pm 1}, \ldots, t_d^{\pm 1}],$$
that is, $c_{\bf m} \neq 0$ for all vertices ${\bf m} \in P$. 
The zero locus 
$$Z := \{ f({\bf t}) = 0 \} \subset 
\T^d_K:={\Spec}\, K[t_1^{\pm 1}, \ldots, t_d^{\pm 1}]$$
we call {\bf affine toric hypersurface}. 

\begin{rem} Note that 
the affine 
toric hypersurface $Z := \{ f({\bf t}) =0 \}\subset \T^d$ determines its defining  non-constant  Laurent polynomial
$f$ uniquely  up to multiplication by a nonzero monomial  
$a {\bf t}^{\bf m}(a \neq 0)$ which shifts the Newton polytope of $f$ by lattice vector ${\bf m} \in M$. So it will be  convenient 
to refer to $P$ as  
{\bf Newton 
polytope} of the toric hypersurface 
$Z \subset \T^d$. Moreover, 
it is natural to consider Newton polytopes $P$ of toric  hypersurfaces $Z \subset \T^d$  up unimodular equivalence, since the group of affine linear transformations 
${\rm Aff}(\Z^d) = GL(d,\Z)  \rtimes  \Z^d$ 
acts on Laurent polynomials $f$ via automorphisms  
${\rm Aut}(\T^d) \cong GL(n,\Z)$ of the algebraic torus 
$\T^d$, and via  multiplication 
by monomials ${\bf t}^{\bf m} = t_1^{m_1} \ldots t_d^{m_d}$.
An unimodular isomorphism   
 $\varphi \in {\rm Aff}(\Z^d) = GL(d,\Z)  \rtimes  \Z^d$ transforms  an affine 
 toric  hypersurface  
$Z_1 \subset \T^d$ with the Newton polytope $P_1$ into the isomorphic affine hypersurface $Z_2 \subset \T^d$ with the Newton polytope $P_2 = \varphi(P_1)$. 
\end{rem}

Now let us consider the geometric meaning of lattice projections $\pi\, : \, P \to P'$ from view point of toric hypersurfaces $Z \subset \T^d$ with the Newton polytope $P$.

\begin{rem}
Let  $P$ be the Newton polytope of a Laurent polynomial $f$. Assume that $P$ admits a standard lattice projection onto a $k$-dimensional lattice polytope $P' \subset \R^{k}$ $(0< k < d)$.  Then we can  view $P'$ as Newton polytope of the Laurent polynomial 
$f $ considered  as element of the 
 Laurent polynomial ring  
 $R[t_1^{\pm}, \ldots, t_{k}^{\pm}]$ whose coefficients  ring is another Laurent polynomial ring 
$R:= \C[t_{k+1}^{\pm 1}, \ldots, t_d^{\pm 1}]$. 
Using the splitting $\T^d \cong  \T^{k} \times \T^{d-k}$ and 
the ring embedding $R \hookrightarrow R [t_1^{\pm}, \ldots, t_{k}^{\pm}]$,  
we obtain the surjective morphism $\T^d \to \T^{d-k} = {\rm Spec}(R)$ whose restriction to $Z$ is a dominant 
 morphism 
$Z \to \T^{d-k}$ such that general fibers are affine toric hypersurfaces in $k$-dimensional torus $\T^{k}$ having $P'$ as   
Newton polytope.  
\end{rem}

Consider some examples of lattice projections.

\begin{exam} \label{01}
Assume that a $d$-dimensional lattice polytope $P$ has {\bf width 1}, that is, $P$ has a lattice projection on the unique
hollow (also $F$-hollow) lattice segment $[0,1] \subset \R$.  The lattice projection $\pi\, :\, P \to P' = [0,1]$
means that  the $d$-dimensional 
lattice polytope $P \subset \R^d$   
is unimodular equivalent to a lattice polytope 
in $\R^d$ contained between two 
parallel integral affine hyperplanes $\{ x_d =0\}$ and 
$\{ x_d =1\}$. Up  to this unimodulal isomorphism, we obtain the corresponding 
Laurent polynomial $f(t_1, \ldots, t_d) \in K[t_1^{\pm 1}, \ldots, t_d^{\pm 1}]$ in the form:
\[ f(t_1, \ldots, t_d) = g_0(t_1, \ldots, t_{d-1}) + t_d g_1(t_1, \ldots, t_{d-1}), \]
for some Laurent polynomials $g_0, g_1 \in K[t_1^{\pm 1}, \ldots, 
t_{d-1}^{\pm 1}]$.  
Since the polynomial $f({\bf t})$ defining $Z \subset \T^d$ is linear with respect to the last variable $t_d$,  we can rationally  eliminate $t_d$ from this equation $f =0$  by the formula 
\[ t_d = -\frac{g_0(t_1, \ldots, t_{d-1})}{g_1(t_1, \ldots, t_{d-1})} \]
and obtain a birational isomorphism 
$Z \stackrel{bir}{\sim} \A^{d-1}_K$ over $K$, i.e., $Z$ 
is an irreducible  $K$-rational algebraic variety.    
\end{exam}

The following conjecture proposes a natural  "inverse statement" to  last example.

\begin{conj} \label{rationality}
Let $P \subset \R^d$ be a $d$-dimensional lattice polytope. 
Assume that  for any field $K$ any toric  hypersurface $Z \subset (K^*)^d$ 
with the Newton polytope $P$ is irreducible and birational to  $\A^{d-1}_K$ 
over $K$. Then the Newton polytope 
$P$ admits a lattice projection onto $[0,1]$.  
\end{conj}

\begin{exam} \label{conic1}
In case $d =2$, the conjecture can be easily verified. Indeed, the rationality of general curve $Z \subset \T^2$ with Newton polygone $P$ implies that $P$ has no interior
lattice points, i.e.,  $P$ is hollow. Up to unimodular isomorphisms, the unique sporadic hollow lattice polygon is   
the  triangle 
\[  Q := {\rm Conv}((0,0), (2,0), (0,2)) \subset \R^2. \]
which is the Newton polytope of a general conic
 $C \subset \A^2_K$ defined by 
a quadratic equation 
\[  a_{0,0} + a_{1,0} t_1  + a_{0,1} t_2 + 
a_{2,0} t_1^2 + a_{1,1}t_1 t_2 + 
a_{0,2} t_2^2 = 0, \;\; (a_{i,j} \in K) . \] 
If $K = \C$, then $C$ 
is birational to $\A^1_\C$. However, we can take $K = \R$ consider  the conic 
\[ 1 + t_1^2 + t_2^2 =0 \]
is not birational to $\A^1_\R$. Moreover, we can consider the conic 
$(1 + t_1)^2 - t_2^2=0$ over any field $K$  with the Newton polygone $T$ consisting  of two irreducible components. 

More generally,  
if a $d$-dimensional  Newton polytope $P$ hypersurface
$Z \subset \T^d$ has a lattice projection onto $Q$, 
then the toric 
hypersurface $Z$ becomes birational to a conic bundle over  
$(d-2)$-dimensional algebraic torus. Note that the corresponding toric hypersurface $Z \subset \T^d$ might 
be non-rational variety even over the algebraically closed field $\C$ (see Example \ref{Klein} below). 
\end{exam}

\begin{exam} \label{Klein}
Let ${\mathcal K}_{d} \subset \P^{d+1}$ be smooth projective $d$-dimensional 
Klein cubic given by  the homogeneous equation 
\[ z_0 z_1^2 + z_1z_2^2 + \cdots + z_d z_{d+1}^2 + 
z_{d+1} z_0^2 = 0, \;\; d \geq 2, \]
which is invariant under the cyclic permutation of the homogeneous coordinates $z_0, z_1, \ldots, z_{d+2}$. 
The Newton polytope $P_{d+1}$ of ${\mathcal K}_d$  
is a weakly sporadic $(d+1)$-dimensional $F$-hollow 
lattice simplex with the minimal multiplier 
$\mu = (d+2)/3$ having no other lattice points besides its vertices.   
The  affine equation of ${\mathcal K}_d \cap \A^{d+1} \subset \A^{d+1} = \{ z_0 =1 \}$:
\[ t_1^2 + t_1 t_2^2 + t_2 t_3^2 + t_3 t_4^2 + \cdots +  t_{d+1} = 0 \]
considered with respect to the pair of variables $x_1$ and $x_3$   defines a conic $C$ from Example \ref{conic1}. 
Note that the lattice $(d+1)$-simplex 
$P_{d+1}$ admits several different lattice projections onto $2$-dimensional simplex $Q:={\rm Conv}((0,0), (2,0), (0,2))$
defining  corresponding 
different conic bundle structures on the Klein cubic   ${\mathcal K}_{d}$. If the dimension $d=2k$ is even, then  $P_{2k+1}$ admits a lattice projection onto $[0,1]$ 
which sends $k+1$ lattice vertices of $P_{2k+1}$ corresponding to cubic monomials $z_{2i}z_{2i+1}^2$ $(0 \leq i \leq k)$ to the lattice point $0$ and 
the remaining $k+1$ lattice vertices of $P_{2k+1}$ to the lattice point $1$.  By \ref{01}, we see that  any even-dimensional Klein 
cubic ${\mathcal K}_{2k}$ is rational 
over any field $K$. Note that it is rather 
nontrivial to show that the  $3$-dimensional 
Klein cubic ${\mathcal K}_3$ is 
not rational over $\C$ (see \cite{CG72}). Moreover, a still open conjecture claims that 
the Klein cubic ${\mathcal K}_{2k+1}$ 
is not rational over $\C$ in all odd dimensions 
$d=2k+1 \geq 3$ (see \cite{Kol19}). This conjecture is supported by the fact that the $2k$-dimensional empty 
lattice simplex $P_{2k}$ has no a lattice projection onto $[0,1]$ for all  integers   
$k \geq 2$.  We note that any $d$-dimensional toric  hypersurface 
$Z \subset \T^{d+1}$ with the Newton polytope $P_{d +1}$ is
irreducible and smooth for $d \geq 3$.  
\end{exam}

\section{Non-degenerate toric hypersurfaces and $\Q$-Fano fibrations}

\begin{dfn}
Let $Z \subset \T^d_K$ be an affine toric hypersurface given as zero locus of a Laurent polynomial 
$$f ({\bf t}) = \sum_{{\bf m} \in A} {c}_{\bf m} {\bf t}^{\bf m} \in 
K[t_1^{\pm 1}, \ldots, t_d^{\pm 1}]$$
with the Newton polytope $P={\rm Conv}(A)$. The toric  hypersurface $Z$ is called {\bf non-degenerate} if $Z$ is smooth over the algebraic closure $\overline{K}$, and for any face $\Theta \prec P$ $(0 < \dim \Theta < d)$ the
affine toric 
hypersurface $Z_\Theta \subset \T^d$ defined as zero locus of the Laurent polynomial 
\[ f_{\Theta}({\bf t}) := \sum_{{\bf m} \in A \cap \Theta} {c}_{\bf m} {\bf t}^{\bf m} \]
is reduced and smooth over $\overline{K}$ as well. 
\end{dfn}

\begin{exam} Let  $P$ be a $d$-dimensional lattice simplex. 
Denote by $A$ the set of its $d+1$ vertices. Assume that ${\rm char}\,K =0$. Then the non-degeneracy of $Z \subset \T^d$ given by $\sum_{{\bf m} \in A} c_{\bf m} {\bf t}^{\bf m}$ is  equivalent to the nonvanishing condition 
\[ \prod_{{\bf m} \in A}  c_{\bf m} \neq 0. \]
\end{exam}

Now we want to explain why Theorem \ref{proj-F} is inspired by the Minimal Model Program for non-degenerate toric hypersurfaces $Z \subset \T^d$ with Newton polytope $P$. 

It was proved in \cite{Ba23} that if the Fine interior $F(P)$ of a $d$-dimensional lattice polytope $P$  is not 
empty, then every non-degenerate toric hypersurface $Z  \subset \T^d$ with 
the Newton polytope $P$ has a minimal model, i.e., a projective model $\widehat{Z}$ with at worst $\Q$-factorial terminal singularities and semi-ample 
canonical $\Q$-divisor.
Moreover, the Fine interior $F(P)$ 
allows to compute  the Kodaira dimension $\kappa(\widehat{Z})$ by the  formula: 
\[ \kappa(\widehat{Z}) = \min \{ d-1, \dim F(P) \}. \]

Now we are interested in birational geometry of non-degenerate toric hypersurfaces $Z \subset \T^d$ with $F$-hollow Newton polytope $P$. 

Our main result is  the following.

\begin{theorem} \label{qFano}
Let $P$ be a $d$-dimensional weakly sporadic  $F$-hollow 
lattice polytope $P$ with the minimal multiplier $\mu =\mu(P)> 1$. Denote by $\P_{Q^*}$ the canonical toric $\Q$-Fano variety whose 
defining fan is spanned by faces of the canonical Fano  polytope
\[ Q:= {\rm Conv}(S_F(\mu P)). \]
We  denote by $\widetilde{Z}$ the Zariski closures in $\P_{Q^*}$
of a non-degenerate toric hypersurface $Z \subset \T^d$ with the Newton polytope $P$. Then $\widetilde{Z}$ is a  
normal variety and the following adjunction formula holds: 
\[ K_{\widetilde{Z}} =  (K_{\P_{Q^*} }+ \widetilde{Z})\vert_{\widetilde{Z}} =\left(\frac{\mu-1}{\mu}\right)K_{\P_{Q^*}} \vert_{\widetilde{Z}}.\] 
 Furthermore, $\widetilde{Z}$ has at worst canonical singularities and ample anticanonical $\Q$-divisor, i.e., 
$\widetilde{Z}$ is a $(d-1)$-dimensional 
canonical $\Q$-Fano hypersurface in  $\P_{Q^*}$.
 \end{theorem}

In the proof we need the following technical statement. 

\begin{lemma} \label{codim2}
Let $P$ be as in \ref{qFano}. Denote by $L(P) \subset \C[t_1^{\pm 1}, \ldots, t_d^{\pm 1}]$ the  linear system on the canonical toric $\Q$-Fano variety 
$\P_{Q^*}$ spanned by all monomials ${\bf t}^{\bf m}$ corresponding to  lattice points  ${\bf m} \in P \subseteq Q^*\cap  M$. Then no codimension-$2$ torus orbit in $\P_{Q^*}$ is contained 
in the base locus of $L(P)$. 
\end{lemma}

\noindent
{\em Proof of Lemma \ref{codim2}.} Take an arbitrary $1$-dimensional edge $E \prec Q$
of the canonical Fano polytope $Q \subset N_
R$. 
Let $\sigma_E \subset N_\R$ be the  
$2$-dimensional cone over $E$ in the fan $\Sigma_{Q^*}$ defining the canonical 
toric $\Q$-Fano variety $\P_{Q^*}$. Since the lattice polytope $Q$ has only one interior lattice point $0 \in N$, the lattice polygon ${\rm Conv}(0, E) \subset Q$ has no lattice points other than $0$ and $E \cap N$. This means  that $\sigma_E$ is unimodular equivalent to  the $2$-dimensional cone $\sigma_k \subset \R^2$ spanned by lattice vectors $(1,0)$ and $(1,k)$ for some $k \geq 1$. Let   $\{e_0, e_1, \ldots, e_k\}:=  E \cap S_F(\mu P)$ be the lattice vectors corresponding under the above equivalence to lattice vectors 
\[ (1,0), (1,1), \ldots, (1,k) \in \sigma_k. \]
Denote by $N' \subset N$ the sublattice in $N$ spanned by 
$e_0, e_1, \ldots, e_k$. Then $N'$ is a direct summand of $N$
since each pair $e_i, e_{i+1}$ forms a part of $\Z$-basis of $N$ for every $i \in \{0, 1, \ldots, k-1 \}$. Take a $\Z$-basis $v_1, v_2, \ldots, v_d$ of $N$ such that $v_1 = e_0$, $v_2 = e_1$ and consider the standard  lattice projection  
$$\pi_E\, :\, M \to M':=Hom(N', \Z) \cong \Z^2$$
together the identification of lattice vectors $e_i \in N'$ with $(1,i) \in \Z^2$ for all $i \in \{0, 1, \ldots, k \}$. Denote by 
 $P'_E \subset \R^2$  the lattice polygon 
$\pi_E(P)$.  Using a shift 
by an appropriate lattice vector ${\bf m} \in M$, 
we can assume without loss of generality that ${\rm Min}_P(e_0) = {\rm Min}_P(e_1) =0$, that is, both hyperplanes $\langle x, e_0 \rangle =0$ and $\langle x, e_1 \rangle =0$ are integral supporting hyperplanes for 
the lattice polytope $P$ and for the rational polytope 
$\mu P$. Since $F(\mu P)$ is just a rational point $p \in M_\Q$,  we obtain 
$$\langle p, e_0 \rangle = {\rm Min}_{F(\mu P)}(e_0)  =  
{\rm Min}_{F(\mu P)}(e_1) = \langle p, e_1 \rangle  =1, $$ 
and by linearity of $\langle p, * \rangle$ we obtain 
$$\langle p, e_i \rangle = {\rm Min}_{F(\mu P)}(e_i) =1, \;\;\; \forall i \in \{0, 1, \ldots, k\}.$$
Since all lattice points $e_0, e_1, \ldots, e_k \in E \prec Q$ are 
contained in $S_F(\mu P)$, we obtain    
 $$ {\rm Min}_{\mu P}(e_i) =0 \;\; \forall i \in \{0, 1, \ldots, k\},$$
i.e., the linear equations $\langle x, e_i \rangle = 0$ define integral supporting hyperplanes for $\mu P$ and for $P$.  

In order to show that the codimension-$2$ torus
orbit corresponding to $\sigma_E$ is not contained in the 
base locus of $L(P)$, we have to show that the $(d-2)$-dimensional linear subspace 
$l_E \subset M_\R$,  
$$l_E := \langle x, e_i \rangle =  0, \;\;\forall i \in \{0, 1, \ldots, k\}, $$  contains at least one lattice vertex of $P$, or equivalently, 
the lattice $\pi_E$-projection $P_E' = \pi_E(P) \subset \R^2$ contains the origin  $(0,0) =\pi_E(l_E)  \in \R^2$. We consider two cases. 

\noindent
{\bf Case 1:} $k \geq 2$. Then $\pi_E$-projection of $P$ is 
a lattice polygon $P_E' \subset \R^2$
contained 
in the cone 
\[ C_k:= \{ (x_1, x_2) \in \R^2 \, \mid \, x_1 \geq 0, x_1 + kx_2 \geq 0 \}. \]
The supporting integral hyperplane  
$\langle x, e_1 \rangle =0$ for $P$ must contain at least one 
vertex of $P$. On the other hand, $\pi_E$-projection of this 
hyperplane is the line in $\R^2$ with the equation $x_1 + x_2 = 0$, and it has  only the origin $(0,0)$ as common point with the cone $C_k$. 

\noindent
{\bf Case 2:} $k=1$.  Then $2$-dimensional cone $\sigma_1 \subset \R^2$ is regular, i.e., its generators $e_0, e_1 \in S_F(\mu P)$ form a $\Z$-basis of $N'$. 
Assume that $(0,0) \not\in P_E' = \pi_E(P)$. Then one of two generators $(0,1)$ and $(1,-1)$ of $\sigma_1$ must 
belong to $P_E'$, since otherwise the lattice polytope 
$P_E'$ would be contained in the convex set 
\[ B:= {\rm Conv}( \{\sigma_1 \cap \Z^2 \}\setminus \{ (0,0), (0,1), (1,-1)\} )  \]
and the lattice point  $(1,0) =\pi_E(F(\mu P))$ would not be in the interior of $\mu P'_E = \pi_E(\mu P) \subset \mu B$ for $\mu >1$. 
So we obtain ${\rm Min}_P(e_0 + e_1) =1$ and ${\rm Min}_{\mu P}(e_0+e_1) = \mu >1$. However, 
we have ${\rm Min}_{F(\mu P)}(e_0 + e_1) = {\rm Min}_{F(\mu P)}(e_0) +  {\rm Min}_{F(\mu P)}(e_1)  =1+ 1 =2$. This implies 
\[ {\rm Min}_{F(\mu P)}(e_1 + e_2) - {\rm Min}_{\mu P}(e_1+e_2)  = 2- \mu <1. \]
The latter contradicts the definition of Fine interior  $F(\mu P)$.  
\medskip

\noindent
{\em Proof of Theorem \ref{qFano}.} 
 Let $\P^{(2)}_{Q^*} \subset \P_{Q^*}$ be the Zariski open toric subvariety obtained   from $\P_{Q^*}$ by deleting all torus orbits of codimension at least $3$ in $\P_{Q^*}$. The Gorenstein 
toric variety $\P^{(2)}_{Q^*}$ has a minimal crepant toric desingularization $\rho\,  :\, \widehat{\P}^{(2)}_{Q^*} \to \P^{(2)}_{Q^*}$ by resolving of possible $A_{k-1}$-singularities along codimension-$2$ strata in 
$\P^{(2)}_{Q^*}$.  Applying  Lemma \ref{codim2}, we obtain  that  the Zariski closure $\widehat{Z}^{(2)}$ of the non-degenerate toric hypersurface $Z$ in $\widehat{\P}^{(2)}_{Q^*}$ 
is smooth, and the Zariski closure  $\widetilde{Z}^{(2)}$ of $Z$ in $\P^{(2)}_{Q^*}$ is Gorenstein, and it has at worst   $A_{k-1}$-singularities 
along the transversal intersections of the non-degenerate quasi-projective toric hypersurface $\widetilde{Z}^{(2)} \subset \P^{(2)}_{Q^*}$ with codimension-$2$ torus orbits in $\P^{(2)}_{Q^*}$. 
Thus we obtain that the singular locus of $(d-1)$-dimensional projective hypersurfaces $\widetilde{Z} \subset \P_{Q^*}$ has codimension at least  $2$ and, by Serre's criterion for  normality,  $\widetilde{Z}$ is normal.  Now we can relate the canonical classes of $\widetilde{Z}$ and $\P_{Q^*}$ by 
the adjunction formula on the Gorenstein quasi-projective 
toric variety $\P^{(2)}_{Q^*}$.

Let $\dim F(\mu P) =0$ and $p= F(\mu P)$. 
We consider the shifted rational polytope $P_0:= \mu P - p$ 
such that $F(P_0) = \{ 0 \}$. Let $\Sigma_{Q^*}[1]$ be the set of generators
of $1$-dimensional cones in the fan defining the canonical 
toric $\Q$-Fano variety $\P_{Q^*}$ . Then $\Sigma_{Q^*}[1] \subseteq S_F(\mu P)$ and we have 
$ {\rm Min}_{P_0}(\nu)  =-1$ for all $\nu \in \Sigma_{Q^*}[1]$. 

The lattice polytope $P$ defines on $\P_{Q^*}$  a toric Weil divisor 
\[ L_P:= \sum_{\nu \in   \Sigma_{Q^*}[1]} (- {\rm Min}_P(\nu))D_{\nu} \]
such that $\Q$-Cartier divisor $\mu L_P$ is rationally equivalent to  
the anticanonical class  $-K_{\P_{Q^*}}$ of $\P_{Q^*}$,
since  
\[ K_{\P_{Q^*}} = \sum_{\nu \in   \Sigma_{Q^*}[1]} 
- D_{\nu}  =  \sum_{\nu \in   \Sigma_{Q^*}[1]}  {\rm Min}_{P_0}(\nu) 
D_{\nu}  =  \sum_{\nu \in   \Sigma_{Q^*}[1]} \left( 
\mu {\rm Min}_{P}(\nu) - \langle p, \nu \rangle \right)  
D_{\nu},      \]
and for any $p \in M_\Q$  the $\Q$-divisor  
\[  \sum_{\nu \in   \Sigma_{Q^*}[1]}  \langle p, \nu \rangle  
D_{\nu}      \]
is a principal $\Q$-divisor of $\P_{Q^*}$. We write the adjunction formula on the Gorenstein toric variety $\P^{(2)}_{Q^*}$ in the form 
\[ K_{\widetilde{Z}} = ( K_{\P_{Q^*}} + L_P)  \vert_{\widetilde{Z}}  = \left( -{\mu} L_P + L_P \right)  \vert_{\widetilde{Z}}  =\left( \frac{\mu -1}{\mu}\right)  K_{\P_{Q^*}} \vert_{\widetilde{Z}}.   \]
In particular, the canonical divisor of $\widetilde{Z}$ is a 
$\Q$-Cartier divisor. It order to show that singularities of 
$\widetilde{Z}$ are at worst canonical, we use a result of Khovanski\^i  
that a non-degenerate toric hypersurface with Newton polytope $P$ always admits smooth projective birational model $W$ in a smooth toric variety defined 
by a regular  refinement $\Sigma$ of the normal fan $\Sigma_P$  \cite{Kho83}. We can always choose a regular refinement $\Sigma$ of $\Sigma_P$ which is simultaneously a regular refinement of $\Sigma_{Q^*}$, so that  
we obtain for smooth toric hypersurface $W \subset \P_\Sigma$ with a birational morphism $\rho\, :\, W \to \widetilde{Z}$ together with  a formula 
\[ K_{W} = \rho^*K_{\widetilde{Z}} + \sum_{\nu \in \Sigma[1]} a_\nu (D_\nu \cap W), \]
where $a_\nu=0$ if $\nu \in S_F(\mu P)$,  or if $D_\nu \cap W = \emptyset$, and $a_\nu = (-1 - {\rm Min}_{P_0}(\nu)) >0$ 
for all $\nu \not\in S_F(\mu P)$. Thus $\widetilde{Z}$ has at worst canonical singularities.  
\hfill $\Box$

\begin{coro} \label{fano-fiber}
Let $P$ be a $d$-dimensional $F$-hollow lattice polytope and $\dim F(\mu P) = k \geq 1$. Denote by 
$\pi\, :\, P \to P'$ the lattice projection onto $(d-k)$-dimensional weakly sporadic $F$-hollow lattice polytope 
$P'$ with $\mu(P') = \mu(P) = \mu$ constructed 
in \ref{proj-F}. Then $\pi$ defines
a dominant morphism $\varphi\,:\,  Z \to \T^k$ whose 
general fibers are  non-degenerate hypersurfaces with the $(n-k)$-dimensional weakly sporadic $F$-hollow lattice polytope $P'$ admitting $\Q$-Fano 
projective compactifications.  
\end{coro}  

\noindent
{\em Proof.} We use Theorem \ref{qFano} in relative situation by considering the $(n-k)$-dimensional 
fan $\Sigma' \subset N'_{\R}$, where $N' \subset N$
is the $(n-k)$-dimensional sublattice from Case 1 in the proof of Theorem \ref{proj-F}. The fan $\Sigma'$ is spanned by 
faces of $(n-k)$-dimensional canonical Fano polytope 
$Q'$, where $Q' = {\rm Conv}(S_F(\mu P)) \cap N')$.
We can identify the finite set $S_F(\mu P)) \cap N'$ with 
the support $S_F(\mu P')$  of the  $(n-k)$-dimensional weakly sporadic $F$-hollow lattice polytope $P' = \pi (P) 
\subset M_\R$ with the minimal multiplier $\mu = 
\mu(P') = \mu(P)$. By Theorem \ref{qFano}, each general 
fiber of the dominant morphism $\varphi\, Z \to \T^k$ admits 
a natural $\Q$-Fano projective compactification which can 
be obtained in the following way. 
  
Using \cite[Theorem 6.3]{Ba23}, 
we embed the non-degenerate toric hypersurface 
$Z \subset \T^d$ into $d$-dimensional $\Q$-Gorenstein 
toric variety $\widetilde{\P}$ associated with the Minkowski sum $F(\mu P) + C(\mu P)$, where $C(\mu P)$ is the $d$-dimensional rational polytope containing $\mu P$  defined 
by the inequalities:
\[ C(\mu P):= \{ x \in M_\R\, \mid \, \langle x, \nu \rangle \geq 
{\rm Min}_{\mu P}(\nu), \;\; \forall \nu \in S_F(\mu P)\}. \]
Since $F(\mu P)$ is a Minkowski summand of $\widetilde{P}$, 
the normal 
fan $\Sigma_{\widetilde{P}}$ contains the $(n-k)$-dimensional fan $\Sigma'$ as a subfan describing 
toric $\Q$-Fano 
fibers of the toric morphism $\alpha \, :\,  \widetilde{\P} \to \P_{F(\mu P)}$, where $\P_{F(\mu P)}$ is 
$k$-dimensional toric variety corresponding to $k$-dimensional rational polytope $F(\mu P)$. By restricting $\alpha$ to the Zariski closure $\widetilde{Z}$ of $Z$ in $\widetilde{\P}$, we obtain the $\Q$-Fano fibration 
$\widetilde{\varphi}\, :\, \widetilde{Z} \to \P_{F(\mu P)}$ which extends the dominant morphism $\varphi\,:\,  Z \to \T^k \subset \P_{F(\mu P)}$ to the hypersurface 
$\widetilde{Z} \subset \widetilde{\P}$.
\hfill $\Box$

\begin{coro} \label{crit}
A non-degenerate hypersurface $Z \subset \T^d$ with  Newton polytope $P$ has negative Kodaira dimension if and only if $P$ is $F$-hollow. 
\end{coro}

\noindent
{\em Proof.} If $F(P) \neq \emptyset$, then $Z \subset \T^d$ has a minimal model which can be constructed as Zariski closure $\widehat{Z} $ in some simplicial torus embedding 
$\T^d \subset \widehat{V}$. In particular, the Kodaira dimension of $Z$ is non-negative \cite{Ba23}.

If $F(P) = \emptyset$, then, by Theorems \ref{proj-F} and \ref{fano-comb},  if $P$ is either weakly sporadic and $Z \subset \T^d$ is birational to a $\Q$-Fano hypersurface in the canonical $\Q$-Fano variety $\P_{Q^*}$, or, by \ref{fano-fiber}, $Z$ is birational to a $\Q$-Fano fibration $\widetilde{Z} \to \P_{F(\mu P)}$ over a $k$-dimensional toric variety. Therefore the projective toric hypersurface $\widetilde{Z}$ has negative Kodaira dimension. 
\hfill  $\Box$

\end{document}